\documentclass{article}
\usepackage[utf8]{inputenc}
\usepackage{amsmath,amsthm,amsfonts,amssymb,amscd}

\usepackage{comment}

\usepackage{dsfont}
\title{{\Large \bf{On the Existence of Word-representable Line Graphs of Non-word-representable Graphs}}}

\author{
M M Akbar\footnote{E-mail: akbar@utdallas.edu},
\,\,
P D Akrobotu\footnote{E-mail: Prosper.Akrobotu@utdallas.edu}
\ \,\,\&\,
C P Brewer\footnote{E-mail: Charlie.Brewer@utdallas.edu}
\\
\\ {{Department of Mathematical Sciences}}
\\ {University of Texas at Dallas}
\\ {Richardson, TX 75080, USA}}

\date{\today}

\newcommand{\mb}[1]{\boldsymbol{#1}}

\usepackage{tikz}
\usepackage{braket}
\usepackage{polynom}
\usepackage{float}
\usepackage[caption = false]{subfig}
\usepackage{ifthen}
\usepackage{graphicx}
\usepackage{hyperref}
\usepackage{esvect}
\newtheorem{thm}{Theorem}

\newtheorem{conjecture}{Conjecture}

\newtheorem{proposition}{Proposition}
\newtheorem{example}{Example}
\begin{document}

\maketitle
\begin{abstract}
An open question in the theory of word-representable graphs for the past decade has been whether the line graph of a non-word-representable graph is always non-word-representable. By formulating an appropriate optimization problem for the decision problem of 3-semi-transitive graphs, we show that the line graph of a non-word-representable graph can be word-representable. 
 Using IBM's CPLEX solver, we demonstrate for several known word-representable and non-word-representable graphs that the line graph of a graph is 3-semi-transitive when there is a solution to the optimization problem. This results in an example where the line graph of a non-word-representable graph is both 3-semi-transitive and semi-transitive and thus is word-representable.
\end{abstract}

\section{Introduction}
 A graph $G=(V,E)$ is word-representable if there exists a word $w$ over the set of alphabets $V$ (set of vertices of $G$) such that letters $x$ and $y$ alternate in $w$ if and only if $\{x,y\}$ is an edge in $E$ (the set of edges of $G$) \cite{KitaevLozin, KP2008}. The word $w$ is then a word-representant of $G$. By letters $x$ and $y$ alternating in $w$, we mean that if we delete all other letters, $w$ reduces to a word of the form $xyxy\ldots xy$ or $yxyx\ldots yx$.
For example, the cycle graph on 5 vertices labeled sequentially by 1, 2, 3, 4 and 5 can be represented by the word 1521324354. Here, the letter $2$ and $3$ alternate (2323) whereas, for example, letters $2$ and $4$ do not (2244). It follows from the definition that for a complete graph $K_{n}$  any permutation of $\{1,2,3,\ldots,\,n\}$ is a word-representant and for an empty graph $E_{n}$, $w=123\cdots (n-1)nn(n-1)\cdots 321$ is a word-representant.

The study of word-representable graphs is relevant to several fields including combinatorics on words and scheduling \cite{KitaevLozin}. Writing down the words of word-representable graphs gets complicated as the number of nodes and edges in the graph gets large \cite{alternation}. Therefore alternative tools like orientation of edges are employed.
An acyclic orientation of edges such that any directed path $v_{1} \rightarrow v_{2}\rightarrow v_{3} \rightarrow \cdots \rightarrow v_{k}$ for $k\geq 4$ in which either the arc $v_{1}\rightarrow v_{k}$ is absent or for every $1\leq i<j\leq k$, the arc $v_{i}\rightarrow v_{j}$ exists is called a semi-transitive orientation. Graphs that admit semi-transitive orientations are called semi-transitively orientable. If any of the conditions above for semi-transitivity is violated, the orientation yields what we call a shortcut graph, i.e., a digraph on the set of vertices $V=\{v_{1},v_{2},\ldots v_{k}\}$ that contains a directed path $v_{1} \rightarrow v_{2}\rightarrow v_{3} \rightarrow \cdots \rightarrow v_{k}$ and an arc $v_{1}\rightarrow v_{k}$ but missing at least the arc $v_{i}\rightarrow v_{j}$ for some $1< i<j< k$ \cite{Kitaev2017}.  Using semi-transitive orientations, it can be shown that a graph is word-representable if and only if it is semi-transitively orientable \cite{HALLDORSSON2016164}. Thus, deciding if a graph with $m$ edges is word-representable requires searching through the space of all $2^{m}$ possible orientations for a semi-transitive orientation. The problem grows exponentially with the number of edges. It has been show that such ``recognition problem'' for word-representable graphs is NP-complete, i.e., a problem with no known efficient procedure for finding solutions \cite{HALLDORSSON2016164}. However, satisfiability modulo theories (SMT) were recently employed in constructing an algorithm to determine the word-representant of a word-representable graph 
in which the authors also introduced ``$k$-semi-transitive orientation'' as a refinement to the notion of semi-transitive orientation to aid in identifying and classifying word-representable graphs \cite{SMT2019}. 
Defining an undirected graph as $k$-semi-transitive, if it admits an acyclic orientation that avoids shortcuts of length $k$ (that is shortcuts of length $\geq k+1$ are allowed), the authors were able to investigate and show using SMT that a $3$-semi-transitive graph on at most $8$ nodes is word-representable. This would be important for this paper. 

The organization of this paper is as follows. In section 2, we present a quadratic constrained binary optimization for detecting 3-semi-transitive line graph of a graph and use IBM CPLEX solver to solve it. In section 3, we verify that the optimization problem correctly recognizes 3-semi-transitive line graphs for a number of graphs and present an explicit example of a word-representable line graph of a non-word-representable graph of maximum degree 4.
We conclude with two conjecture.

\section{Optimization problem}
The formulation of the decision problem (of a word-representable graph) that we present below can be solved using both classical solvers and quantum algorithms. In short, using semi-transitive orientation or $3$-semi-transitive orientation we will construct a satisfiability problem, i.e., a problem that requires searching for a state $\mb{x}=(x_{1},\ldots,x_{n})$ of $n$ variables that satisfy a conjunction $C(\mb{x}) = \bigwedge_{i=1}^{m}C_{i}$ of $m$ clauses $C_{i}$ (Boolean functions) on $n$-variables $x_{i}$, as an optimization problem.
We want to determine whether there is a given instance (set of clauses) with an assignment that violates no clauses.

Let $G=(V,E)$ be an undirected graph with $V(G) = \{1,2,\ldots,n\}$ and $E(G) = \{e_{1},e_{2},\ldots,e_{m}\}$ as the set of $n$ vertices and $m$ edges of the graphs, respectively.
Define a vector $\mb{x} = (x_{1},x_{2},\ldots,x_{m})$ of dimension $m\times 1$ with entries $x_{i} \in \{-1,1\}$. 
Let $M$ be the $0$--$1$ incidence matrix of $G$. Then we consider its edge adjacency matrix $Q = M^{T}M$ for $G$ and formulate the word-representability of the line graph $L(G)$ of $G$ as a minimization problem whose constraints are forced to capture the properties of  $3$-semi-transitive orientations. 

\begin{proposition}
\label{claim1}
Let $G$ be an undirected graph with incidence matrix $M$ such that the edge adjacency matrix of $G$ is $Q=M^{T}M -2I$ where $I$ is the identity matrix. Then the line graph $L(G)$ with adjacency matrix $Q$ is $3$-semi-transitively orientable if the minimization problem equation (\ref{eq:spin_min_prob_const}) is solvable, where $\mb{x}=(x_1,\dots,x_m)$ and $m$ is the number of edges of $G$.
\begin{align}
\label{eq:spin_min_prob_const}
    \min_{\mb{x}}\, \hspace{5pt} &\mb{x}^{T}Q\mb{x}\\
    s.t.\,
    &C_{1}:\,  |x_{i} + x_{j} + x_{k}| = 1 \,\text{ if } q_{ij}= q_{ik}=q_{jk} = 1 \text{ for } i\neq j\neq k\nonumber\\
    &C_{2} \,: \, x_{i}\,\in\,\{-1,1\}\nonumber
\end{align}
\end{proposition}

The optimization problem can be converted into a binary optimization problem using the transformation $\mb{x} = 2\mb{y} - \mb{1}$  where $\mb{y}$ is a tuple of binary variables $y_{i}\in \{0,1\}$ and $\mb{1}$ is a vector of 1's. Therefore the binary equivalent of equation (\ref{eq:spin_min_prob_const}) is given by
 \begin{align}
\label{eq:min_prob_const}
    \min_{\mb{y}}\, &4\mb{y}^{T}Q\mb{y} + \mb{b}^{T}\mb{y} + c \\
     s.t.\,
    &C_{1}:\,  1\leq y_{i} + y_{j} + y_{k} \leq 2 \,\text{ if } q_{ij}=q_{ik}=q_{jk} = 1 \text{ for } i\neq j\neq k\nonumber\\
    &C_{2} \,: \, y_{i}\,\in\,\{0,1\}\nonumber
\end{align}
where $ \mb{b} = 4(Q-2I)\mb{1}$ and $c = \mb{1}^{T}\mb{b}.$

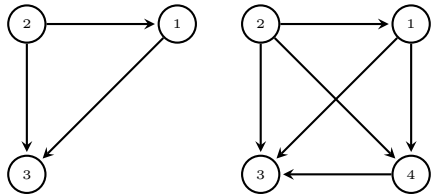
\begin{figure}[!h]
    \centering
    \begin{tikzpicture}[->,>=stealth,shorten >=1pt,auto,node distance=2cm,
                        thick,main node/.style={circle,draw,font=\sffamily\tiny\bfseries}]
    
     \node[main node] (1) {$1$};
      \node[main node] (2) [left of=1] {$2$};
      \node[main node] (3) [below  of=2] {$3$};
      
      \draw[every node/.style={font=\sffamily\small}]
       
        (2) edge node[right] {} (1)
        (2) edge node[left] {} (3)
        (1) edge node[right] {} (3);
        
    \end{tikzpicture}
    \hspace{10pt}
      \begin{tikzpicture}[->,>=stealth,shorten >=1pt,auto,node distance=2cm,
                        thick,main node/.style={circle,draw,font=\sffamily\tiny\bfseries}]
    
     \node[main node] (1) {$1$};
      \node[main node] (2) [left of=1] {$2$};
      \node[main node] (3) [below  of=2] {$3$};
      \node[main node] (4) [below of=1] {$4$};
      
      \draw[every node/.style={font=\sffamily\small}]
       
        (4) edge node[left] {} (3)
        (2) edge node[left] {} (3)
        (1) edge node[left] {} (3)
        (2) edge node[left] {} (1)
        (1) edge node[left] {} (4)
        (2) edge node[left] {} (4);
    \end{tikzpicture}
    \caption{Subgraph captured by constraint $C_{1}$ (left) and a 3-semi-transitive graph (right).}
    \label{fig:figG2}
    \end{figure}

To proceed with verifying this proposition, we first recall that the line graph of an undirected graph $G$ is a graph $L(G)$ describing the adjacency between the edges of $G$. That is, the adjacency matrix $Q$ of $L(G)$ has a non-zero entry $q_{ij}=1$ if and only if the edges $e_{i}$ and $e_{j}$ are incident to a common node and $q_{ij}$ is zero otherwise. Thus, the nodes of $L(G)$ are the edges of $G$. Secondly, from the definition of $k$-semi-transitive graphs we have that an undirected graph is $3$-semi-transitive if it admits an acyclic orientation such that for any directed path
$v_{0} \rightarrow v_{1} \rightarrow v_{2} \rightarrow v_{3}$ of length 3, if $v_{0} \rightarrow v_{3} $ is an edge, then so are $v_{0} \rightarrow v_{2}$ and $v_{1} \rightarrow v_{3}$.
The constraint $C_{1}$ captures subgraphs of the line graph $L(G)$ that are isomorphic to the complete graph $K_{3}$ on 3 vertices (triangles) and assigns weights $x_{i}\, \in \,\{-1,1\}$ to the nodes of the subgraph where $-1$ signifies a sink and $1$ signifies a source such that there is always one pure sink and one pure source between every 3 adjacent nodes just as show in Figure \ref{fig:figG2}. Here, node $2$ is a pure source and node $3$ is a pure sink, however nodes $1$ and $4$ are mixed. This implies a typical solution of the optimization problem will have $x_{1}=x_{2}=1$ and $x_{3}=x_{4}=-1$. This means that any given solution would, in the worst case, generate a partial orientation of the edges $\{i,j\}\in E(L(G))$ in the line graph $L(G)$ of the graph $G$ for all $i,j$ such that $x_{i}\neq x_{j}$. Therefore to assign an orientation to the edges $\{i,j\}\in E(L(G))$ with $x_{i}=x_{j}$ to obtain a 3-semi-transitive orientation, we employ the Boolean function defined in \cite{SMT2019} : \begin{equation}
\label{eq:sergey3-semi}
        \left(\exists\, k,\,m \,:\, (e_{ik}\land e_{jk})\lor(e_{ki}\land e_{kj})\land e_{im} \land e_{mj} \right)\implies e_{ij},
    \end{equation}
    where $e_{ik}$ is a Boolean variable set to true if there is a directed edge from $i$ to $k$ and false otherwise. That is, if the condition is satisfied then a 3-semi-transitive orientation requires an edge directed from $i$ to $j$. Using this definition and the acyclic orientation condition, a 3-semi-transitive orientation of the line graph can be determined from a solution of the optimization problem. 

Clearly, if the line graph $L(G)$ of $G$ is 2-colorable, then we will have an ideal solution giving rise to an alternating pure source and pure sink orientations yielding a semi-transitive orientation. However if $L(G)$ is 3-colorable then the possible orientations obtained from a solution are isomorphic to those shown in Figure \ref{fig:figG2}.
\begin{example}\rm
Consider the complete graph $G=K_{4}$ shown in Figure \ref{fig:example} (a), with line graph $L(G)$ shown in Figure \ref{fig:example} (b). We seek to find a semi-transitive orientation using the method outlined in this paper. Solving the constrained optimization problem equation (\ref{eq:spin_min_prob_const}) with CPLEX, we find the following solution for $\mb{x}$:
\begin{equation}
    \mb{x}=(-1,1,1,-1,-1,-1).
\end{equation}
where
\begin{equation}
    V(L(G))=\left((0, 1), (1, 2), (0, 3), (2, 3), (0, 2), (1, 3)\right).
\end{equation}
We use this vector to assign each node of the line graph as a source or sink node. Then we assign any edge of $L(G)$ from source to sink, if possible. This gives the following partial orientation shown in Figure \ref{fig:example} (c). However, this still leaves the question of how to orient the edges between adjacent sources or adjacent sinks. To do so, we apply the following algorithm to enforce 3-semi-transitivity:
\begin{itemize}
    \item Choose an edge $(i,j)\in E(L(G))$ which is not yet directed.
    \item Check if either orientation of the edge would cause a directed cycle. If yes, choose the other orientation. If neither orientation would complete a cycle, continue to the next step.
    \item Check if there exist $k,m\in V(L(G))$ which satisfy condition (\ref{eq:sergey3-semi}). If such a pair is found, then the edge is directed from $i$ to $j$.
    \item If no such $k,m$ are found, then edge can be oriented in either direction. Then, move on to the next edge.
\end{itemize}
Applying this algorithm to the current example gives the following orientation in Figure \ref{fig:example} (d) which is easily checked to be 3-semi-transitive. Note that in some cases, the algorithm may result in an orientation which is actually semi-transitive (which is stronger than 3-semi-transitive), however this is not generally guaranteed.

\begin{figure}
    \subfloat[Graph $G =K_{4}$]{\includegraphics[width=0.25\linewidth]{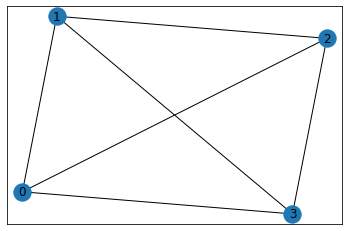}}
    \subfloat[Line graph of $G$]{\includegraphics[width=0.25\linewidth]{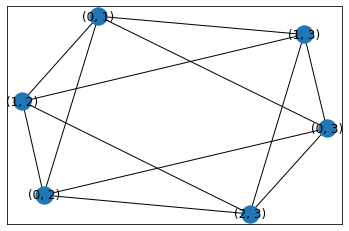}}
    \subfloat[Partial orientation]{\includegraphics[width=0.25\linewidth]{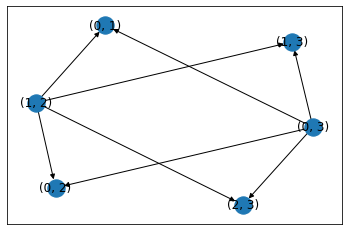}}
    \subfloat[Final orientation]{\includegraphics[width=0.25\linewidth]{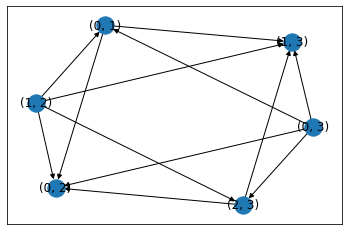}}
    \caption{The complete graph $K_{4}$ and a 3-semi-transitive orientation of its line graph.}
    \label{fig:example}
\end{figure}
\end{example}

\section{Results and Discussions}
\begin{figure}[!h]
\subfloat[Tutte graph $G$]{\includegraphics[width = 1.5in]{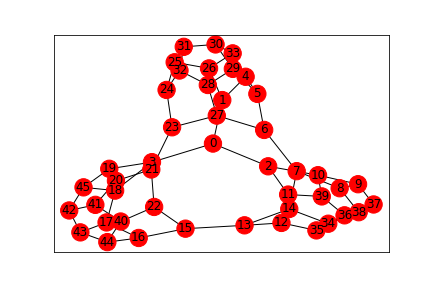}}
\subfloat[Hexagonal lattice]{\includegraphics[width = 1.5in]{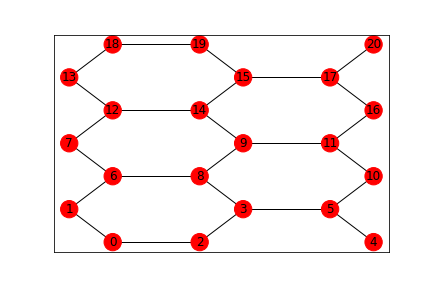}}
\subfloat[Graph $J_4$ from \cite{SMT2019} ]{\includegraphics[width=1.5in]{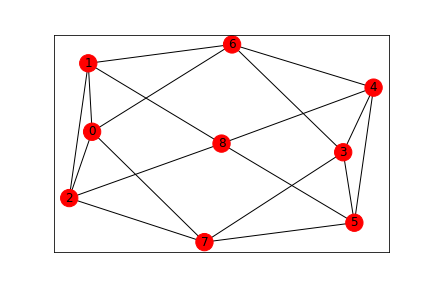}}\\
\subfloat[Wheel graph $W_{4}$]{\includegraphics[width = 1.5in]{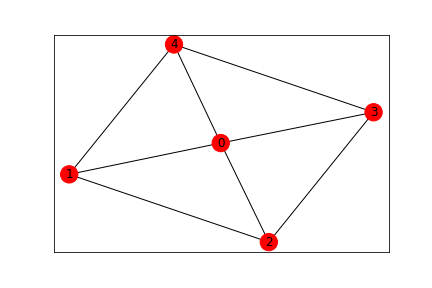}}
\subfloat[Herschel graph ]{\includegraphics[width = 1.5in]{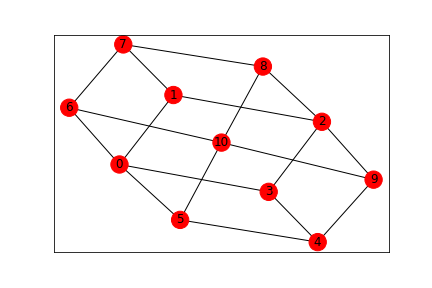}}
\subfloat[Medial graph of Herschel graph $G_{19}$]{\includegraphics[width = 1.5in]{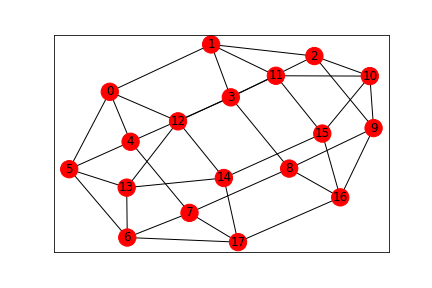}}\\
\subfloat[Complete graph $K_{4}$]{\includegraphics[width = 1.5in]{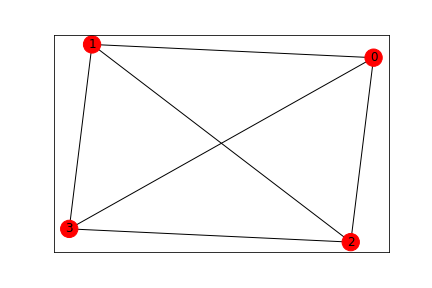}}
\subfloat[Complete graph $K_{4}$ with a broken link]{\includegraphics[width = 1.5in]{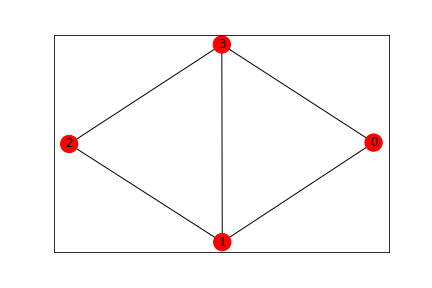}}
\subfloat[Cycle graph $C_{4}$]{\includegraphics[width = 1.5in]{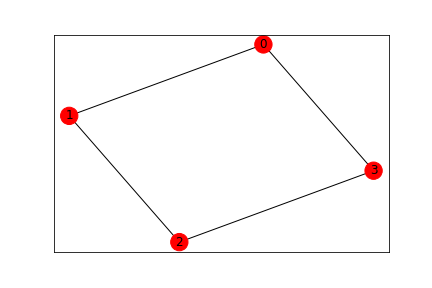}}\\
\subfloat[Petersen graph ]{\includegraphics[width = 1.5in]{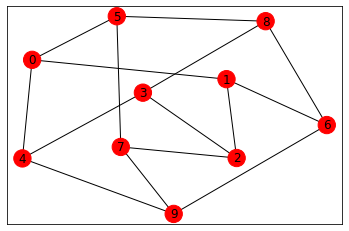}}
\subfloat[Path graph $P_{5}$]{\includegraphics[width = 1.5in]{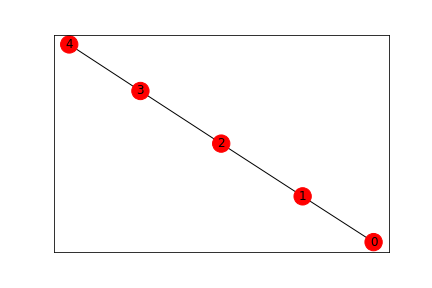}}
\caption{Some semi-transitive graphs whose line graphs are identified to be 3-semi-transitive by proposition \ref{claim1}.}
\label{fig_semi_transitive}
\end{figure}

\begin{figure}[!h]
\subfloat[Graph $T_{1}$ from \cite{pros2015}]{\includegraphics[width = 1.5in]{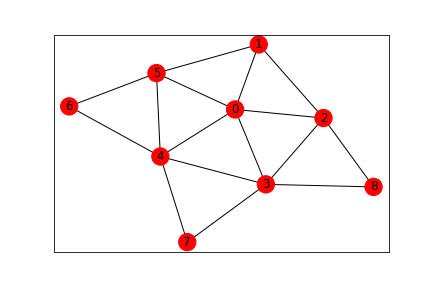}}
\subfloat[Graph $T_{2}$ from \cite{pros2015}]{\includegraphics[width = 1.5in]{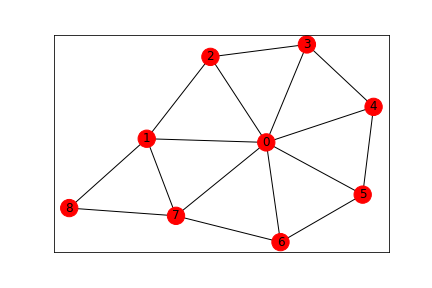}}
\subfloat[Wheel graph $W_{5}$]{\includegraphics[width = 1.5in]{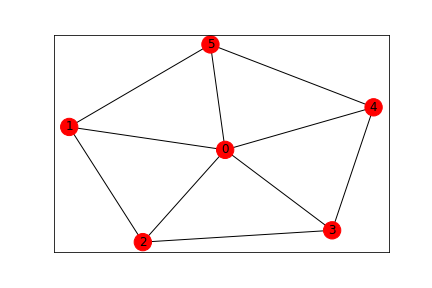}}\\
\subfloat[]{\includegraphics[width = 1.5in]{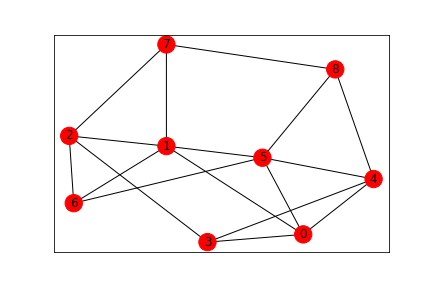}}
\subfloat[]{\includegraphics[width = 1.5in]{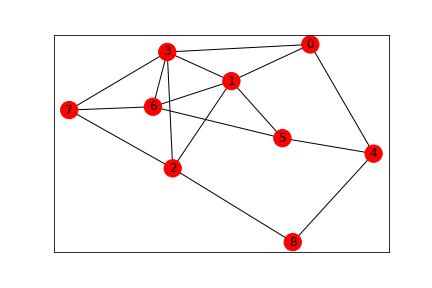}}
\subfloat[]{\includegraphics[width = 1.5in]{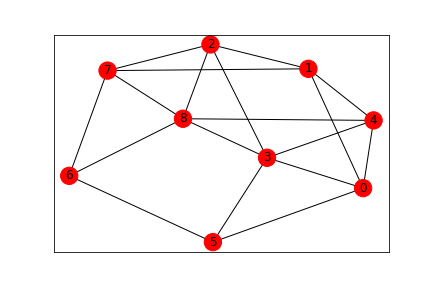}}\\
\subfloat[]{\includegraphics[width = 1.5in]{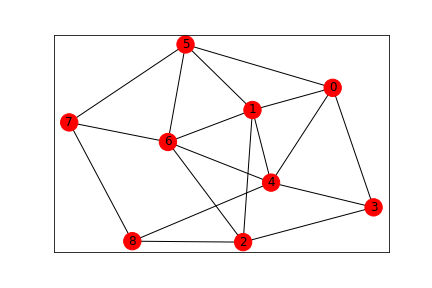}}
\caption{Some non-semi-transitive graphs whose line graphs are not 3-semi-transitive by proposition \ref{claim1}. Graphs from (d) - (g) are from Figure 7 of \cite{SMT2019}.}
\label{fig:non_semi_transitive}
\end{figure}

To verify our proposition, the optimization problem was implemented for several word-representable and non-word-representable graphs (see Figure \ref{fig_semi_transitive}, and Figure \ref{fig:non_semi_transitive} using Python packages \texttt{NumPy}~\cite{numpy}, 
\texttt{NetworkX}~\cite{networkx},  \texttt{Matplotlib}~\cite{matplotlib}, and \texttt{Qiskit}~\cite{qiskit} and 
\texttt{IBM CPLEX} solver \cite{cplex}. We considered several word-representable and non-word-representable graphs and used IBM's  \texttt{CPLEX} tools to generate a quadratic program for the constrained optimization problem given by equation (\ref{eq:spin_min_prob_const}) for each graph and solve using the CPLEX solver.

Among the several graphs considered are polyhedral graphs: 3-vertex connected planar graph such as the Tutte graph, a 3-colorable and 3-regular graph of 46 nodes and 69 edges, the Goldner-Harary graph, a 3-connected maximal planar graph with $11$ vertices, $27$ edges and chromatic number 4, the Hershel graph, a 3-vertex-connected undirected planar graph with $11$ nodes, $18$ edges and chromatic number $2$ and the Medial graph of the Hershel graph.  The formulation was also examined for the hexagonal lattice, path graph, cycle graph, star graph and wheel graphs. We note that all 2- and 3-colorable graphs such as the hexagonal lattice, the Tutte graph, Hershel graph, Medial graph of the Hershel graph, path, cycle and star graphs are all word-representable and equivalently semi-transitively representable graphs. The wheel graphs $W_{2k+1},\, k\geq 2$ with even number of nodes are however not semi-transitively orientable (equivalently non-word-representable) with the smallest non-word-representable graph being the $W_{5}$, the wheel graph with $6$ nodes. It is known that line graphs of the wheel graph $W_{n}$ and complete graphs $K_{n}$ are non-word-representable for $n\geq 5$ \cite{line_graph}.

\begin{table}[!h]
\centering
\resizebox{\textwidth}{!}{
\begin{tabular}{|p{5em}|p{5em}|p{5em}|p{5em}|p{5em}|p{5em}|p{5em}|}
\hline 
\multicolumn{5}{|c|}{Graph Features} & \multicolumn{2}{|c|}{Decision Problem of STO} \\ 
\hline 
Name & \#Nodes & \#Edges & Chromatic No. & Max. deg & STO & QCBO (Solvable) \\
\hline 
Path $P_{n}$ & $n$ & $n-1$ & 2 & 2 & Yes & Yes \\
\hline 
Cycle $C_{n}$ & $n$ & $n$ & $2$ & $2$ & Yes & Yes \\
\hline 
Wheel $W_{2k+1},\,k\geq 1$ & $2k+2$ & $4k+2$ & $4$ & $2k+1$, $k\geq 1$ & No ($k\geq 2$) & No $(k\geq 2)$\\
\hline 
Wheel $W_{2k}, \, k\geq 1$ & $2k+1$ & 4k & 3 & $2k+1$ & Yes & No $(k\geq 3)$ \\
\hline 
Hexagonal lattice ($2\times 3$)  & 21 & 25 & 2 & 3 & Yes & Yes \\
\hline 
Tutte & $46$ & 69 & 3 & 3 & Yes & Yes \\
\hline 
Goldner Harary & 11 & 27 & 4 & 8 & • & No\\
\hline 
Herschel & 11 & 18 & 2 & 4 & Yes & Yes \\
\hline 
Medial of Herschel & 18 & 36 & 4 & 4 & Yes & Yes \\
\hline 
Petersen & 10 & 15 & 3 & 3 & Yes & Yes \\
\hline 
$K_{n}$ & $n$ & $2(n-2)$ & $n$ & $n-1$ & Yes & Yes ($n\leq 5$) \\
\hline 
$J_{4}$ & 9 & 18 & 3 & 4 & Yes & Yes \\
\hline 
$T_1$ & 9 & 16 & 4 & 5 & No & No \\
\hline 
$T_{2}$ & 9 & 16 & 4 & 6 & No & No\\
\hline 
$A$ & 8 & 12 & 4 & 4 & No & Yes\\
\hline
\end{tabular}}
\caption{The decision problem of 3-semi-transitive orientable (3-STO) graphs from quadratic constrained binary optimization (QCBO)}
\label{tabresult}
\end{table}

\begin{figure}[!h]
    \includegraphics[width=.16\textwidth]{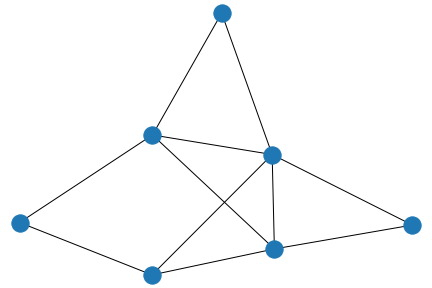}\hfill
    \includegraphics[width=.16\textwidth]{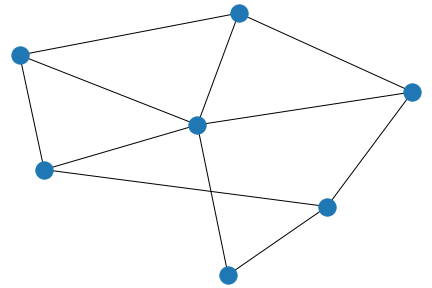}\hfill
    \includegraphics[width=.16\textwidth]{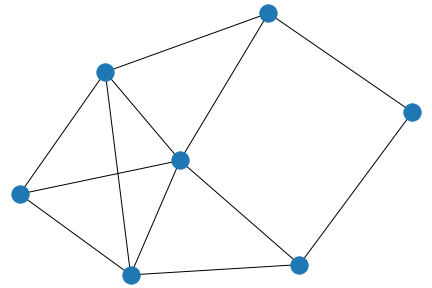}\hfill
    \includegraphics[width=.16\textwidth]{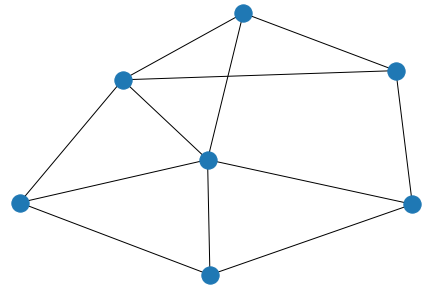}\hfill
    \includegraphics[width=.16\textwidth]{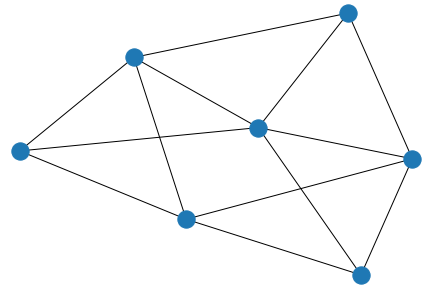}\hfill
    \includegraphics[width=.16\textwidth]{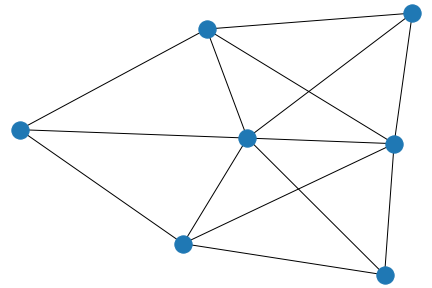}\hfill
    \\[\smallskipamount]
    \includegraphics[width=.16\textwidth]{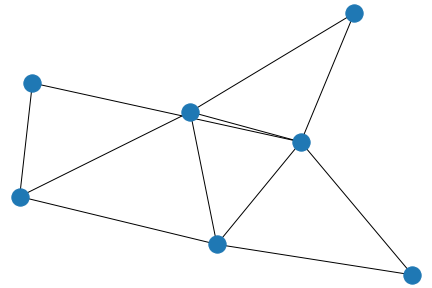}\hfill
    \includegraphics[width=.16\textwidth]{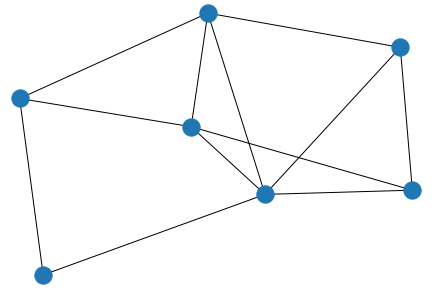}\hfill
    \includegraphics[width=.16\textwidth]{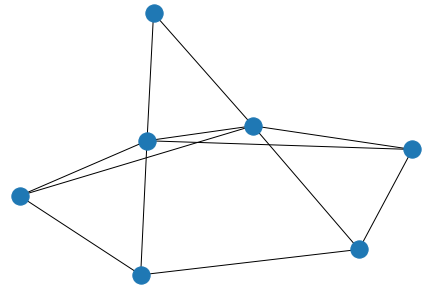}\hfill
    \includegraphics[width=.16\textwidth]{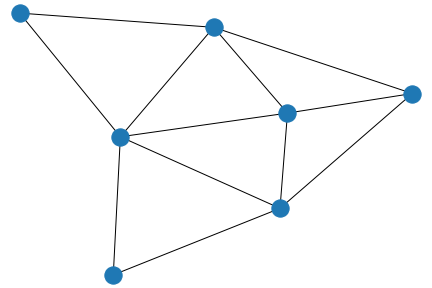}\hfill
    \includegraphics[width=.16\textwidth]{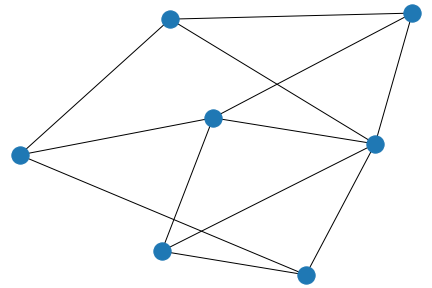}\hfill
    \includegraphics[width=.16\textwidth]{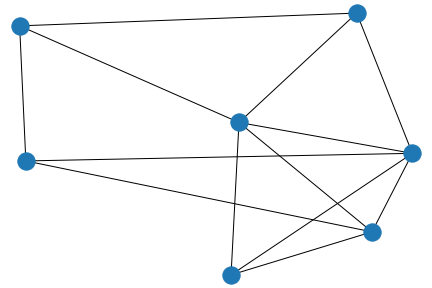}\hfill
    \\[\smallskipamount]
    \includegraphics[width=.16\textwidth]{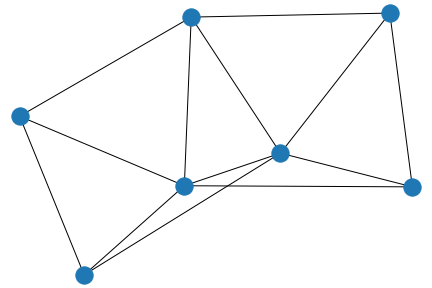}\hfill
    \includegraphics[width=.16\textwidth]{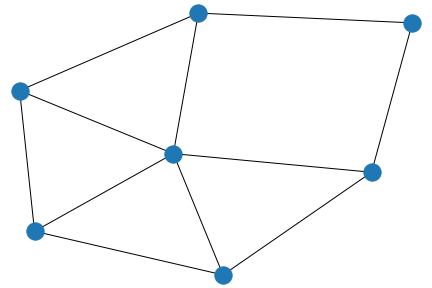}\hfill
    \includegraphics[width=.16\textwidth]{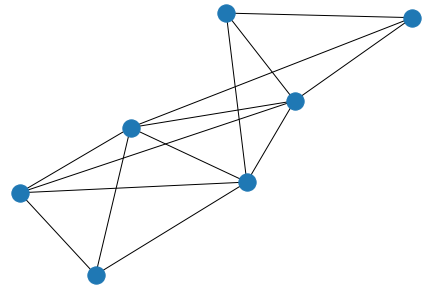}\hfill
    \includegraphics[width=.16\textwidth]{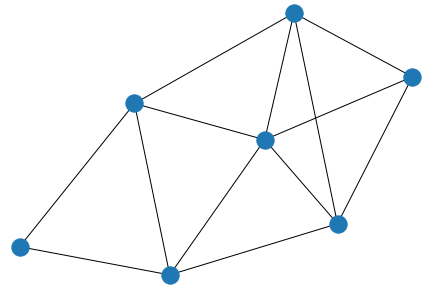}\hfill
    \includegraphics[width=.16\textwidth]{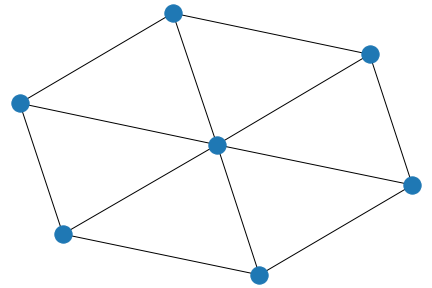}\hfill
    \includegraphics[width=.16\textwidth]{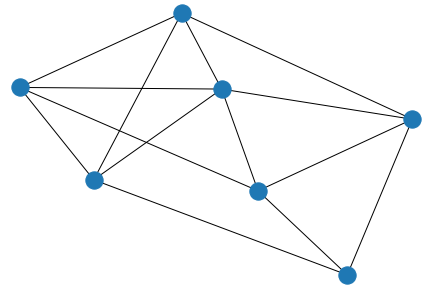}\hfill
    \caption{Some semi-transitive graphs whose line graphs are not 3-semi-transitive by proposition \ref{claim1}.}
    \label{fig:non-3semitrans}\bigskip
    
    \includegraphics[width=.16\textwidth]{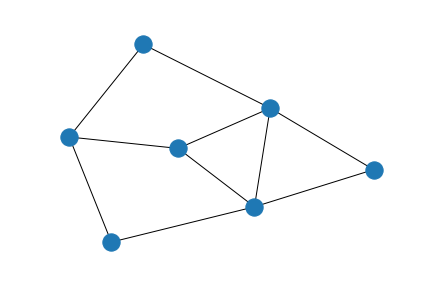}\hfill
    \includegraphics[width=.16\textwidth]{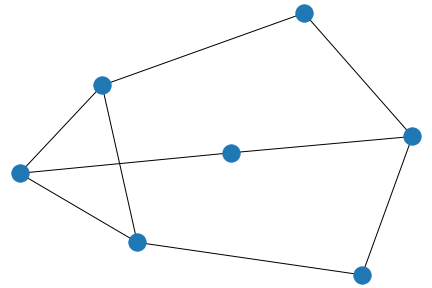}\hfill
    \includegraphics[width=.16\textwidth]{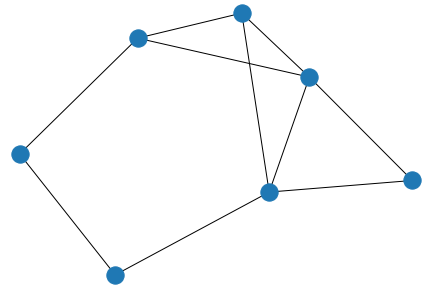}\hfill
    \includegraphics[width=.16\textwidth]{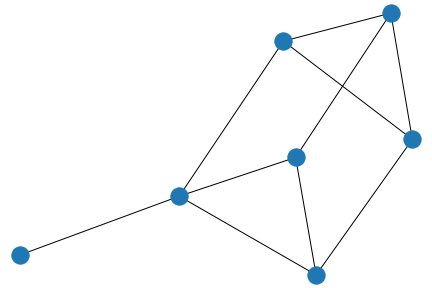}\hfill
    \includegraphics[width=.16\textwidth]{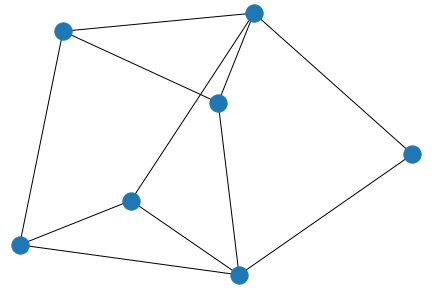}\hfill
    \includegraphics[width=.16\textwidth]{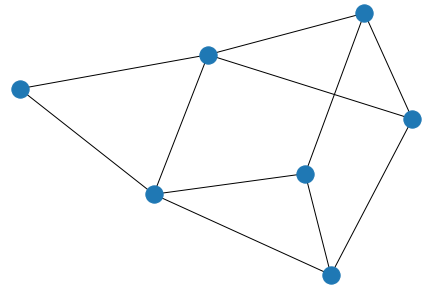}\hfill
    \\[\smallskipamount]
    \includegraphics[width=.16\textwidth]{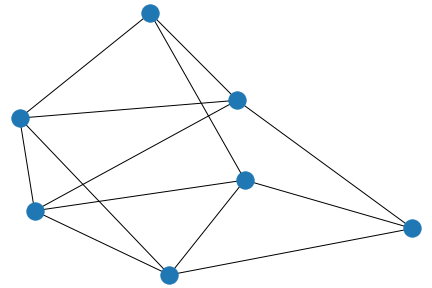}\hfill
    \includegraphics[width=.16\textwidth]{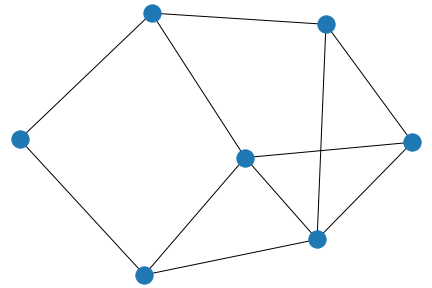}\hfill \includegraphics[width=.16\textwidth]{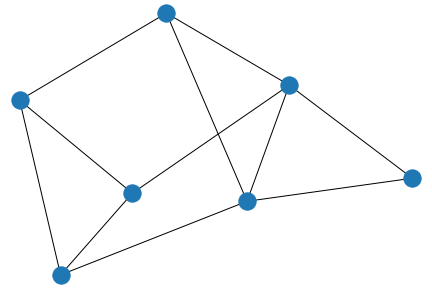}\hfill
    \includegraphics[width=.16\textwidth]{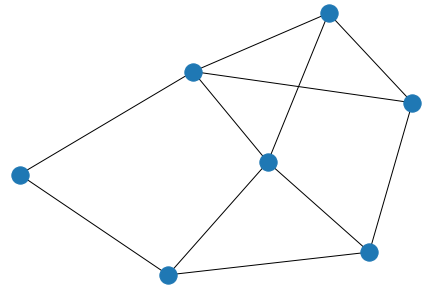}\hfill
    \includegraphics[width=.16\textwidth]{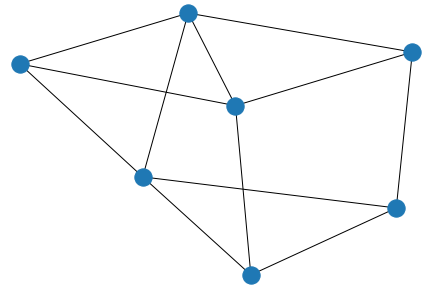}\hfill
    \includegraphics[width=.16\textwidth]{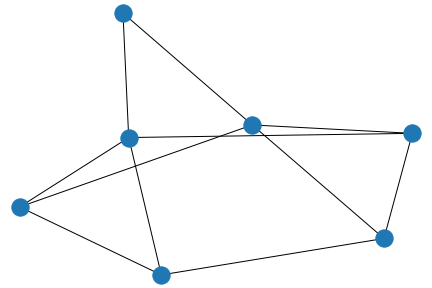}\hfill
    \\[\smallskipamount]
    \includegraphics[width=.16\textwidth]{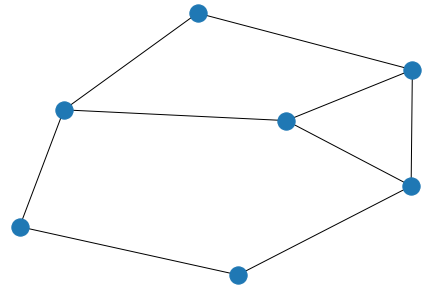}\hfill
    \includegraphics[width=.16\textwidth]{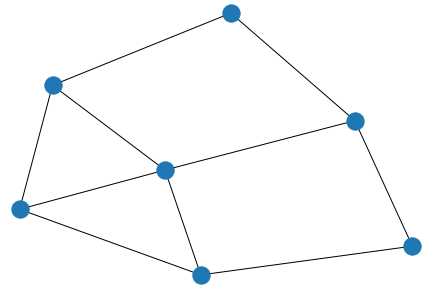}\hfill \includegraphics[width=.16\textwidth]{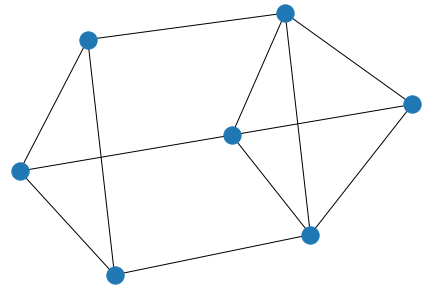}\hfill
    \includegraphics[width=.16\textwidth]{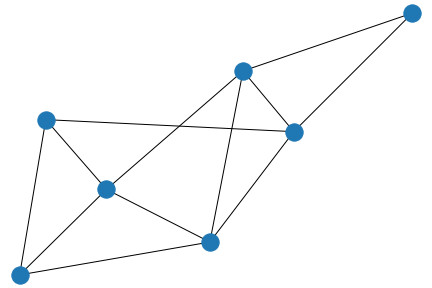}\hfill
    \includegraphics[width=.16\textwidth]{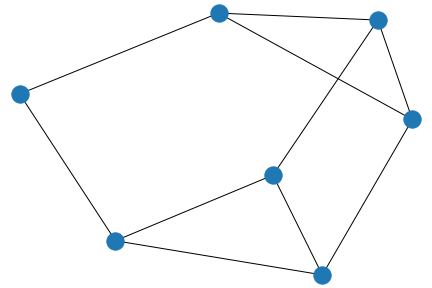}\hfill
    \includegraphics[width=.16\textwidth]{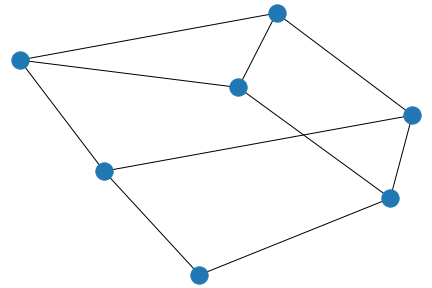}\hfill
    \\[\smallskipamount]
    \includegraphics[width=.16\textwidth]{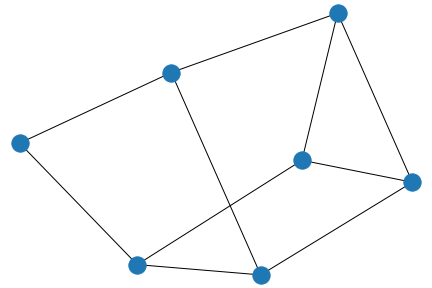}\hfill
    \includegraphics[width=.16\textwidth]{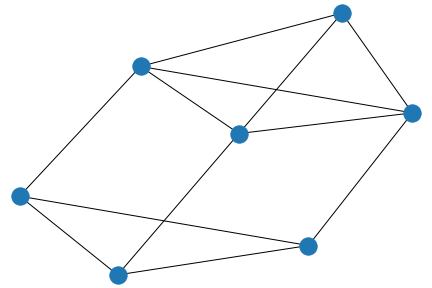}\hfill 
    \hspace{0.5\textwidth}
    \caption{Some semi-transitive graphs whose line graphs are identified to be 3-semi-transitive by proposition \ref{claim1}.}
    \label{fig:3semitrans}
\end{figure}

The NetworkX package was used in creating the graphs considered. Initial attempts were performed on the IBM Q Experience \cite{Q_Experience} platform using the classical solvers in Qiskit. Once the graphs are created, a quadratic program is constructed from the constrained optimization problem of the objective function equation (\ref{eq:min_prob_const}) using Qiskit's CPLEX tools to encode the problem variables, objective function and constraints. The CPLEX solver is then initiated to solve for the minimum solution to the optimization problem. The results obtained using the IBM CPLEX solver showed that the constrained optimization problem was solvable for the semi-transitive graphs in Figure \ref{fig_semi_transitive}  and Figure \ref{fig:3semitrans} 
but not solvable for the following:
\begin{enumerate}
    \item The semi-transitive graphs 
    shown in \cite{SMT2019} and Figure \ref{fig:non-3semitrans}.
    \item The non-word-representable graphs in Figure \ref{fig:non_semi_transitive}. 
    \item Star and wheel graphs of at least $6$ nodes.
\end{enumerate} 
All these graphs have a node of degree at least $5$ which implies the existence of subgraphs isomorphic to the star graph $S_{k}$ of at least $k+1$ nodes. This leads us to the hypothesis that the optimization problem equation (\ref{eq:spin_min_prob_const}) is solvable for a graph $G$ if the maximum degree of $G$ is $\Delta (G)\leq 4$. (This would lead to the conjecture that we present later in the paper.)

\subsubsection*{Word-representable Line Graph of a Non-word-representable Graph}

\begin{figure}[!h]

    \subfloat{\includegraphics[width = 0.25\linewidth]{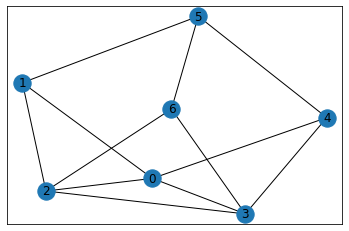}}
    \subfloat{\includegraphics[width=0.25\linewidth]{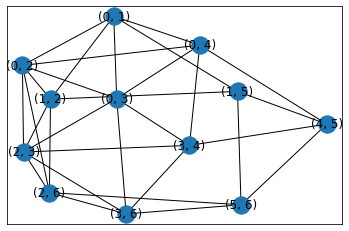}}
    \subfloat{\includegraphics[width=0.25\linewidth]{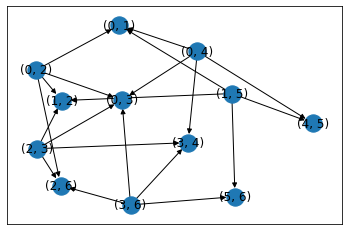}}
    \subfloat{\includegraphics[width=0.25\linewidth]{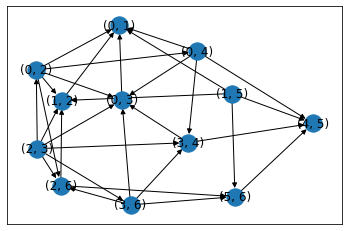}}
    \caption{A non-word-representable graph ($A$) (left) whose line graph (second from the left) is identified to be 3-semi-transitive (right) by proposition \ref{claim1}. The third Figure from the left is the partial orientation obtained from the solution $(1, -1,  1,  1, -1, -1,  1, -1, -1,  1, -1, -1)$ of the optimization problem.}
    \label{fig:graphA}
\end{figure}

The examination of the above hypothesis leads us to the observation that the graph shown in Figure \ref{fig:graphA} is solvable which implies that its line graph is 3-semi-transitive. 

It is easy to check that 3-semi-transitive orientation shown in Figure \ref{fig:graphA} is in fact a semi-transitive orientation and that the line graph of the non-word-representable graph is semi-transitive. 
As mentioned earlier, whether the line graph of non-word-representable graphs are always non-word-representable was an open problem \cite{Kitaev2017} first published in \cite{line_graph} 10 years ago. What we have found here is that the line graph of a graph is word-representable  if the graph is of maximum vertex degree 4. Since the line graph of the graph in Figure \ref{fig:graphA} is solvable and its 3-semi-transitive graph is semi-transitive it answers the open question posed above. It is also important to note that the line graph of wheel graphs on at least $6$ nodes are not word-representable \cite{Kitaev2017} and by our results presented here they are neither 3-semi-transitive. 

\begin{thm}
The line graph $L(A)$ of the non-word-representable graph $A$ is word-representable. 
\end{thm}

\subsubsection*{Conjecture}
Now that we know that $L(G)$ of $G$ can be word-representable while $G$ is not, we propose the following possible relationships between the semi-transitivity of any graph $G$ and its line graph $L(G)$:
\begin{conjecture}
The line graph of a graph $G$ with $\Delta(G)\leq 4$ is at least 3-semi-transitively orientable. The line graph of a non-word-representable graph $G$ is word-representable if $\Delta(G)\leq 4$ and $G$ is planar.
\end{conjecture}
As we mentioned before, line graphs of wheel graphs $W_{n}$ and complete graphs $K_{n}$ are not word-representable for $n\geq 5$ and hence $\Delta (G)\geq 5$ is kept out of the scope of the conjecture \cite{line_graph}. 
\noindent

\section{Conclusion}
The initial objective of this investigation was to obtain a quadratic unconstrained binary optimization (QUBO) formulation for the recognition problem of word-representable graphs (an NP-complete problem) that can be encoded and solved using both classical and quantum computers. 
We considered the sub-problem of recognising  3-semi-transitive graphs and recast the sub-problem as a quadratic constrained binary optimization problem of the line graph $L(G)$ of graph $G$. We found using IBM CPLEX solver that the 3-semi-transitivity of $L(G)$  depends on the solvability of the optimization problem. It is important to note here that the solution to the optimization problem does not directly produce possible 3-semi-transitive orientations of the line graph. However, by devising an algorithm that takes as its input the output of the optimization problem (more precisely, the bit-string or spin configuration of the latter) one obtains a particular 3-semi-transitive orientation of $L(G)$. The explicit example, given in Figure \ref{fig:graphA}, found for the graph $A$ (in \cite{COLLINS2017136}), shows that the line graph on non-word-representable graphs is not always non-word-representable, thus settling the open question known for about a decade which motivated this investigation. 
\section*{Acknowledgement}
We acknowledge the International Business Machines Corporation for providing IBM CPLEX for free which has been used for this work.
\bibliographystyle{ieeetr}
\bibliography{References1}

\begin{thebibliography}{10}

\bibitem{KitaevLozin}
S.~Kitaev and V.~Lozin, {\em Words and Graphs}.
\newblock Springer, Cham, 2015.

\bibitem{KP2008}
S.~Kitaev and A.~Pyatkin, ``On representable graphs,'' {\em Journal of
  Automata, Languages and Combinatorics}, vol.~13, no.~1, pp.~45--54, 2008.

\bibitem{alternation}
M.~M. Halld{\'o}rsson, S.~Kitaev, and A.~Pyatkin, ``Alternation graphs,'' in
  {\em Graph-Theoretic Concepts in Computer Science} (P.~Kolman and
  J.~Kratochv{\'i}l, eds.), (Berlin, Heidelberg), pp.~191--202, Springer Berlin
  Heidelberg, 2011.

\bibitem{Kitaev2017}
S.~Kitaev, ``A comprehensive introduction to the theory of word-representable
  graphs,'' in {\em Developments in Language Theory. DLT 2017. Lecture Notes in
  Computer Science} (R.~M. Charlier~\'{E}., Leroy~J., ed.), vol.~10396,
  pp.~36--67, Springer, Cham, 2017.

\bibitem{HALLDORSSON2016164}
M.~M. Halldórsson, S.~Kitaev, and A.~Pyatkin, ``Semi-transitive orientations
  and word-representable graphs,'' {\em Discrete Applied Mathematics},
  vol.~201, pp.~164--171, 2016.

\bibitem{SMT2019}
{\"O}.~Akg{\"u}n, I.~Gent, S.~Kitaev, and H.~Zantema, ``Solving computational
  problems in the theory of word-representable graphs,'' {\em Journal of
  Integer Sequences}, vol.~22, pp.~1--18, Feb. 2019.

\bibitem{numpy}
``Numpy,'' 2020.
\newblock Online; accessed December 25, 2020.

\bibitem{networkx}
A.~A. Hagberg, D.~A. Schult, and P.~J. Swart, ``Exploring network structure,
  dynamics, and function using networkx,'' in {\em Proceedings of the 7th
  Python in Science Conference (SciPy 2008)}, SciPy 2008, pp.~11--16, ACM,
  2008.

\bibitem{matplotlib}
J.~D. Hunter, ``Matplotlib: A 2d graphics environment,'' {\em Computing in
  Science \& Engineering}, vol.~9, no.~3, pp.~90--95, 2007.

\bibitem{qiskit}
``Qiskit,'' 2020.
\newblock Online; accessed December 25, 2020.

\bibitem{cplex}
``Converters for quadratic programs,'' 2020.

\bibitem{pros2015}
P.~Akrobotu, S.~Kitaev, and Z.~Mas\'{a}rov\'{a}, ``On word-representability of
  polyomino triangulations,'' {\em Sib. Adv. Math.}, p.~1–10, Feb 2015.

\bibitem{line_graph}
S.~Kitaev, P.~Salimov, C.~Severs, and H.~Ulfarsson, ``Word-representability of
  line graphs,'' {\em Open Journal of Discrete Mathematics}, vol.~01, no.~02,
  p.~96–101, 2011.

\bibitem{Q_Experience}
``Ibm quantum experience,'' 2020.

\bibitem{COLLINS2017136}
A.~Collins, S.~Kitaev, and V.~V. Lozin, ``New results on word-representable
  graphs,'' {\em Discrete Applied Mathematics}, vol.~216, pp.~136--141, 2017.
\newblock Special Graph Classes and Algorithms — in Honor of Professor
  Andreas Brandstädt on the Occasion of His 65th Birthday.

\end{thebibliography}
\end{document}